\documentclass[11pt,a4paper]{article}

\usepackage{xcolor}
\definecolor{hellmagenta}{rgb}{1,0.5,0}

\usepackage{color}
\usepackage{amscd}
\usepackage{times}
\usepackage{latexsym}
\usepackage{amsmath}
\usepackage{amssymb}
\usepackage{dsfont}
\usepackage{mathrsfs}
\usepackage{epsfig}
\usepackage{graphicx,graphics}
\usepackage{subfigure}
\usepackage{pstricks}
\usepackage{pst-all}
\usepackage[body={15cm,23cm}, top=2cm]{geometry}
\usepackage{booktabs}

\usepackage{multirow}
\usepackage{diagbox}
\usepackage{url}

\usepackage{bm}
\usepackage{enumerate}
\usepackage{mathrsfs}
\usepackage{amsmath,amssymb,amsthm}

\usepackage{amsmath,esint}
\newtheorem{exmp}{Example}
\usepackage{caption}
\usepackage{cite}
\usepackage{algorithm}
\usepackage{algorithmic}

\usepackage[T1]{fontenc}
\usepackage[utf8]{inputenc}
\usepackage{authblk}

 \newcommand*{\affaddr}[1]{#1} 
 \newcommand*{\affmark}[1][*]{\textsuperscript{#1}}

\begin{document}
	\date{}
	\title{Numerical integration over implicitly defined domains with topological guarantee}
	\author{Tianhui Yang \affmark[1],
		 Ammar Qarariyah \affmark[1], 
		 Hongmei Kang\affmark[2],
		 Jiansong Deng\affmark[1]\thanks{Corresponding author: dengjs@ustc.edu.cn},\\
		\affaddr{\affmark[1]School of Mathematical Science, University of Science and Technology of China, China}\\
	\affaddr{\affmark[2]School of Mathematical Sciences, Soochow University, China}\\
}
\maketitle

%



\section*{Abstract}
Numerical integration over the implicitly defined domains is challenging due to topological variances of implicit functions. In this paper, we use interval arithmetic to identify the boundary of the integration domain exactly, thus getting the correct topology of the domain. Furthermore, a geometry-based local error estimate is explored to guide the hierarchical subdivision and save the computation cost. Numerical experiments are presented to demonstrate the accuracy and the potential of the proposed method.\\
\textbf{Key Words:} Isogeometric Analysis; Numerical integration; Implicitly defined domains; Topological guarantee; Interval arithmetic; Local error estimate; Hierarchical subdivision

\section {Introduction}	
In recent years, the rapid development of 3D printing drives the study of related topics on modeling and analysis. Since implicit representation is convenient in describing the slicing structures of objects \cite{huang2013intersection}, there has been increasing attention paid to representation in this form. Certain properties related to implicit geometries are specifically beneficial for 3D printing. For example, the inside/outside point position test can be applied naturally using an implicit function. Detecting the place of each point in the domain is crucial for the modern layer by layer printing process. Moreover, material features such as coloring can be incorporated implicitly in the object using blending operations which is ideal for biological 3D painting applications (e.g. in~\cite{wegst2015bioinspired}). Building on such capabilities, implicit representation can provide in many cases a 3D printing friendly geometrical models with straightforward constructions.

Integrating CAD and CAE will reduce the huge percent of time cost in the communications between the two processes, thus it has been the goal considered in some work \cite{cottrell2009isogeometric, upreti2017signed, dokken2018trivariate, hollig2003finite, rank2012geometric}. Numerical integration forms an important part of such approaches~\cite{xu2017improved}. Several quadrature techniques have been proposed to ensure an accurate implementation for analysis~\cite{bartovn2017gauss, barendrecht2018efficient, xu2017isogeometric}. For merits of the same idea, the analysis for shapes defined in implicit representation should also be kept in the implicit form.  
This invokes the classical problem of numerical integration in analysis. The only difference is that the domain here is defined implicitly.

When considering the numerical integration over an implicitly defined domain, two problems need to be solved. The first problem is the correct topology of the integration domain due to the topology variance of implicit representation. Generally, finer partitions are used to preserve the correct topology, but it is time-consuming. The second problem is how to get an accurate result of integration. Since many methods consider the approximated mesh as the exact domain to facilitate the calculation, the approximated region leads to significant errors. The boundary condition for analysis through this way is not satisfied exactly.

In this paper, we use the implicit representation itself during the integration process, and try to get the correct topology automatically.
The contributions of the paper are as follows:
\begin{itemize}	
	\item Tackle the topological problem: we employ interval arithmetic in the hierarchical frame to attain the correct topology of the implicitly defined domain.
	\item Save the integration points: a geometry-based local error estimate is proposed to guide the hierarchical subdivision. Such a criterion saves the numbers of integration points when achieving the same accuracy compared with methods not using criterion. Furthermore, the pessimistic error estimate leads to results always smaller than the given tolerance. In this way, we somehow know the "distance" to the exact integral.
	
	\item Theoretical error estimate: a theoretical error estimate is given, where the integration error is caused by the approximation to the boundary of the domain.
\end{itemize}

The organization of the remaining sections is as follows. In Section \ref{re}, we review strategies for integrals on implicitly defined geometry. 
The proposed method is presented in Section \ref{method}. In Section \ref{tests}, several numerical experiments are presented to show the accuracy and efficiency. Concluding remarks and possible future researches are given in Section \ref{con}.

\section{Related work} \label{re}

Numerical integration over domains defined by iso-values of implicit functions have been considered in many methods in science and engineering computing methods. For example, in weighted extended B-splines method \cite{hollig2003finite}, immersed boundary methods \cite{dolbow1999numerical,peskin2002immersed,li2006immersed}, extended finite element methods \cite{moes1999finite,sukumar2001modeling,cheng2010higher}, etc. In this section, we will briefly review the main ideas concerning this topic. For results concerning the numerical integration over rectangular shapes, the readers are refered to \cite{Edalat1999Numerical,berntsen1991adaptive,davis2007methods}.

There are mainly three kinds of ideas to tackle such integration. The first is to use the divergence theorem thus rewriting the original integral as a boundary integral \cite{muller2013highly}. The second idea is to use Monte Carlo method \cite{press2007numerical,shapiro2002architecture}, which can be used for integration over an arbitrary domain. 
Another idea is to modify the integration domain, either by converting the original domain into a rectangular domain enclosing it with the help of discrete Dirac delta function or Heaviside functions \cite{tornberg2004numerical}, or by explicitly reconstructing the boundary\cite{rvachev1994numerical,cheng2010higher,thiagarajan2014adaptively,thiagarajan2017shape}.

Moment fitting method is widely used under the first kind of idea \cite{muller2013highly}, where the numerical coefficients are obtained with a group of given integration points and divergence free functions. In \cite{thiagarajan2014adaptively, thiagarajan2016adaptively}, moment fitting method is performed over a simplified homeomorphic domain, then shape sensitivity analysis is taken to improve the results.  
However, polygonal divergence free function basis must be provided. Also, the quadrature weights may not be positive, thus leading to unstable stiffness matrix in analysis.

Monte Carlo integration method only requires the signs of the function values and avoids the curse of dimensionality, so it can be used for integration over arbitrary domains and higher dimensional spaces \cite{press2007numerical}. However, due to its probability background, it requires large numbers of points and also exhibits slow convergence. Therefore, this method is not appropriate if accurate results are required \cite{shapiro2002architecture}.

Using discrete Dirac delta function or Heaviside function to smear the domain into a larger regular one leads to easier calculations. Since there are many ways to get the integration value over a regular domain. However, it depends on local representations and sometimes does not derive accurate results because of the discontinuous function in the integrand \cite{tornberg2004numerical}.

When explicitly reconstructing the implicit boundary, piecewise linear interpolations are frequently used \cite{rvachev1994numerical, olshanskii2016numerical}. For better approximating the domain, higher order of approximations to the boundary or hierarchical partition are adopted \cite{cheng2010higher}. However, all the methods mentioned above only considered techniques over the simplified polygon. In \cite{engwer2016geometric}, a topology preserving marching cube is proposed, which still considers the discrete level-set function as the exact geometry. In \cite{hollig2015programming}, the integration results are obtained by firstly deriving the intersection points in each direction recursively, then adopting numerical quadrature at the explicit region. Both the shape and the topology are considered in recent work \cite{thiagarajan2018shape}, where specific corrections are added to increase the accuracy, but there are still some limitations or prescribed feature assumptions. That is, the topology could not be well preserved by the initial approximation and later corrections are made. In \cite{saye2015high} and \cite{drescher2017high}, the integration over hyperrectangles and hypersurfaces are considered, respectively.

Here we introduce a bit more details about integration obtained by explicitly reconstruction with hierarchy. The result is the summation of numerical integration over subdivided cells. Cells are classified into three categories: interior, boundary, and exterior cells. Each cell is classified according to the signs of function values at the four corner points. The boundary cells will be recursively partitioned until some criteria are satisfied. At the last level of subdivision, different numerical schemes will be taken according to the type of the cell. The termination criterion is set either by the posterior error  or doing recursive subdivisions until predefined subdivision depth is reached \cite{gautschi1989gauss, press2007numerical}.

In this paper, a numerical integration method based on the hierarchical subdivision is put forward. It tackles the topology problem in implicitly defined domain, by using interval arithmetic. Besides, it uses the implicit representation itself instead of using the discrete mesh approximation to the boundary of the domain. The new termination criterion saves the number of integration points compared with hierarchy-based methods using subdivisions to a predefined level. 

\section{Method}\label{method}
Given an implicit function $f: \mathbb{R}^2\to\mathbb{R}$, which defines a domain $\Omega=\{(x,y)|f(x,y)\geq0\}$, our goal is to keep the correct topology while obtaining the numerical result of $\int_\Omega F(x,y)dxdy$, where $F(x,y)$ is a smooth integrand. During the process, a tolerance $\tau$ is also given to guide the accuracy of the integration.

The hierarchical idea is employed together with two tools. The first is interval arithmetic, which helps to maintain the completeness of integration region, while the second is the subdivision criterion, which attempts to save the number of subdivisions. In this section, we will firstly introduce interval arithmetic, then list the main steps of the proposed method.

\subsection{Interval arithmetic}
Interval arithmetic was introduced by Ramon E.Moore in the 1960s and is widely used in numerical analysis, computer graphics, geometric modeling and global optimization \cite{moore1966interval, mitchell1990robust, moore2014reliability}. One important use of interval arithmetic is in rendering implicit curves and surfaces~\cite{martin2002comparison} with several imporovments on the original algorithim~\cite{shou2003modified}. The idea is to use interval instead of simple numbers to make accurate computations \cite{moore1979methods,gomes2009implicit}.

An interval number $X$ defines a set of numbers: 
$$X=[x_a,x_b]=\{x|x_a\leq{x}\leq{x_b}, x_a\leq x_b\}.$$
	
Given two interval numbers $A=[a,b]$ and $B=[c,d]$, the corresponding rules for the four basic operations are:
\begin{itemize}	
	\item Addition: $A+B=[a+c,b+d]$
	\item Subtraction: $A-B=[a-d,b-c]$
	\item Multiplication: $A\cdot B=[\min\{a\cdot c,a\cdot d,b\cdot c,b\cdot d\},\max\{a\cdot c,a\cdot d,b\cdot c,b\cdot d\}]$
	\item Division: $A/B=[\min\{a/c,a/d,b/c,b/d\},\max\{a/c,a/d,b/c,b/d\}],0\notin[c,d]$
\end{itemize}

In scientific computing, interval quadrature has also been considered, where the remainder of Taylor series is taken on intervals to get the accurate integral \cite{moore1966interval}. 
Different with this usage, we employ interval arithmetic to detect the small feature in implicitly defined domains. As will be seen in the next section, interval arithmetic is used to find enclosures for the ranges of functions in each cell \cite{caprani1981integration,wolfe1998interval}, thus getting its correct sign.

\subsection{Sub-cell classification} \label{class}
Before we explore in more detail, let us first get the correct classifications of the cells at each level of the subdivision.

In hierarchy-based schemes, the type of a cell is usually determined by function values at its four corner vertices. That is, if all the function values are positive (negative), then this cell is an interior (exterior) cell; or else, it is a boundary cell. However, judgments using corner function values are not enough to tell the right types of cells sometimes, especially for cells in implicitly defined domains.

\begin{figure}
	\newcommand{\imgwidth}{3.8in}
	\centerline{
		\includegraphics[width=\imgwidth]{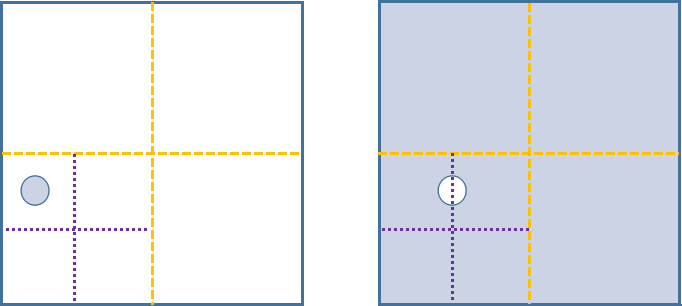}
	}
	\caption{Incorrect classifications of cells: Both cells (with blue frame) should be boundary cells. Only with function values at corner points lead to wrong types: exterior (left) and interior (right), respectively. The interior domains are blue painted.
	}\label{fig:intervalgood}
\end{figure}

Figure \ref{fig:intervalgood} shows two cells with misjudged types. Take the figure on the left as an example. It shows a boundary cell, as we can see that there is a small circle inside. But when just using corner values, we could not get the correct sign even subdividing the cell twice. The wrong classifications will lead to inaccurate integration results. 
In this paper, interval arithmetic is adopted to locate the boundary better.  
 
It proceeds as follows: Assume that the boundary of the domain is defined by $f(x,y)=0$, and $C=[a,b]\times [c,d]$ is a cell whose sign is to be determined. 
Let $X=[a,b]$, $Y=[c,d]$, we calculate the interval number $Z=f(X,Y)$. If $Z$ is a positive (negative) interval, then $C$ is an interior (exterior) cell; otherwise, it is a boundary cell. What important is that, if $0\in Z$, the boundary of the domain is indicated to cross this cell, hence $C$ is a boundary cell. This property guarantees that the boundary will not be missed. 
Now applying this strategy to the domains in Figure \ref{fig:intervalgood}. Both intervals of function values straddle zero, thus we clarify that both cells are boundary cells.

Having derived the correct classifications of the cells, accurate numerical integration could be obtained over interior cells following the tensor product Gaussian quadrature. We now turn our attention to numerical integration over boundary cells.

\subsection{Boundary approximation}\label{appro}
To make the integration over a boundary cell accurate, we shall first get an appropriate representation of the boundary curve. Usually, approximation in explicit form facilitates the repeated numerical integration since it avoids the mass root-finding process in implicit representation.

\begin{figure}[!h]
	\newcommand{\imgwidth}{4.2in}
	\centerline{
		\includegraphics[width=\imgwidth]{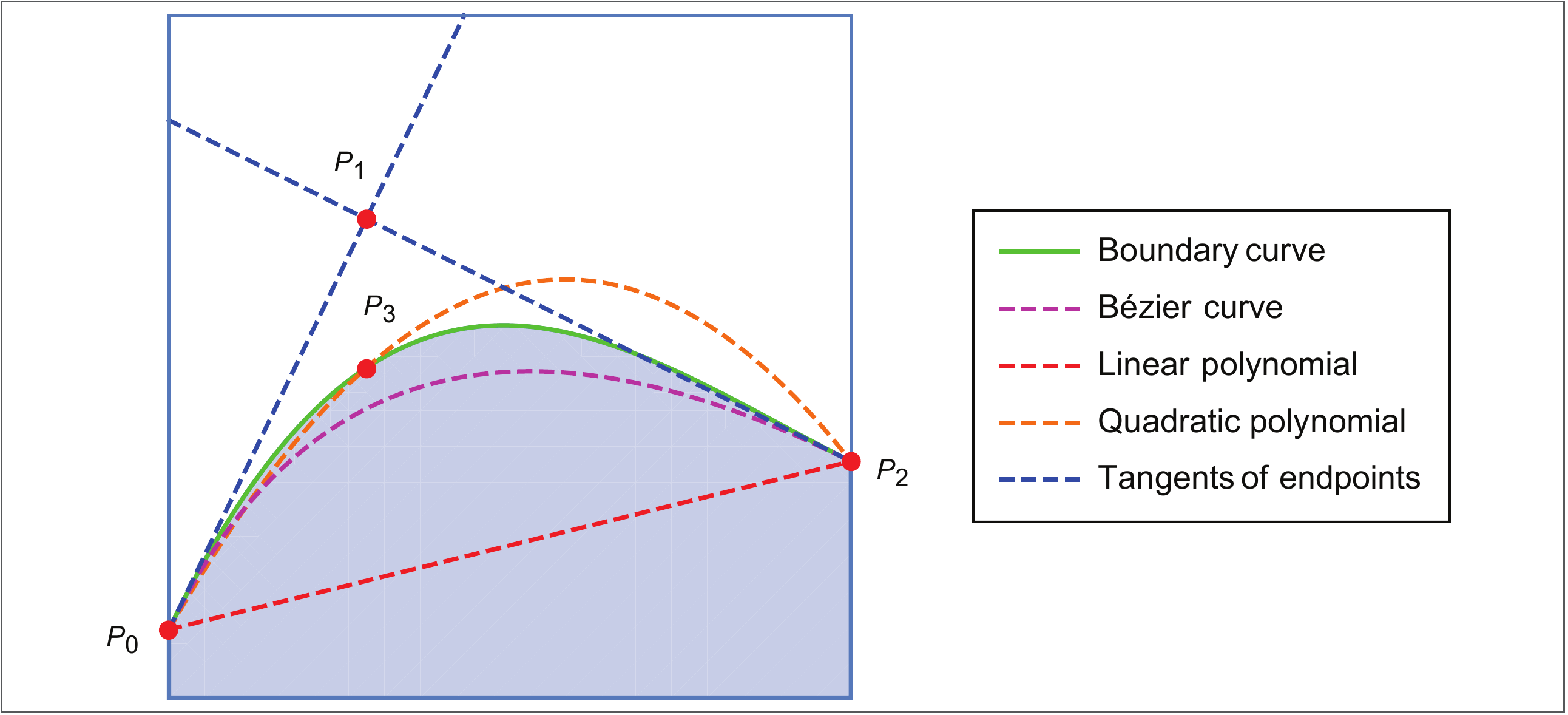}
	}
	\caption{Different approximations to the boundary curve of the implicitly-defined domain, where $P_0$ and $P_2$ are intersections of the boundary curve over the cell, $P_1$ is the intersection of tangents at the two intersections, and $P_3$ is a point on the boundary curve.
	}\label{fig:beziercmp}
\end{figure}

For a clear understanding, three different ways are described when approximating the boundary curve over a boundary cell. As illustrated in Figure \ref{fig:beziercmp}, the simplest way is to connect the two intersections directly, thus deriving the linear polynomial approximation. After using another point $P_3$ together, a quadratic approximation is achieved. Furthermore, a quadratic B\'{e}zier approximation can be determined by using control points $P_i,i=0,1,2$. 

In fact, implicit function theorem indicates the derivative information locally at the two intersections (if $\partial{f}/\partial{y}\neq0$). Moreover, the interpolation property and the derivative relation at the endpoints $P_0$ and $P_2$ convey that $P_1$ is just the intersection point of the tangents at these two points \cite{farin2002curves}. B\'{e}zier form also leads to convenient calculations, as we will see in the next subsection.

Note that the B\'{e}zier approximation could also handle the case if there are singular points at the boundary (if the singular points could be known in advance).
The method is to split the domain at the singular points and getting the tangents at different subregions with the help of blow up\cite{hartshorne2013algebraic}. An example will be taken to show this singular case.

\subsection{Local error estimate guiding further subdivisions over boundary cells} \label{error}
Both linear and quadratic approximation above are usually carried out under a given level of subdivisions, thus leading to excessive computations. 
To avoid this, a geometrically based subdivision rule is proposed to decrease the numbers of calculation.

\begin{figure}[!h]
	\newcommand{\imgwidth}{4.8in}
	\centerline{
		\includegraphics[width=\imgwidth]{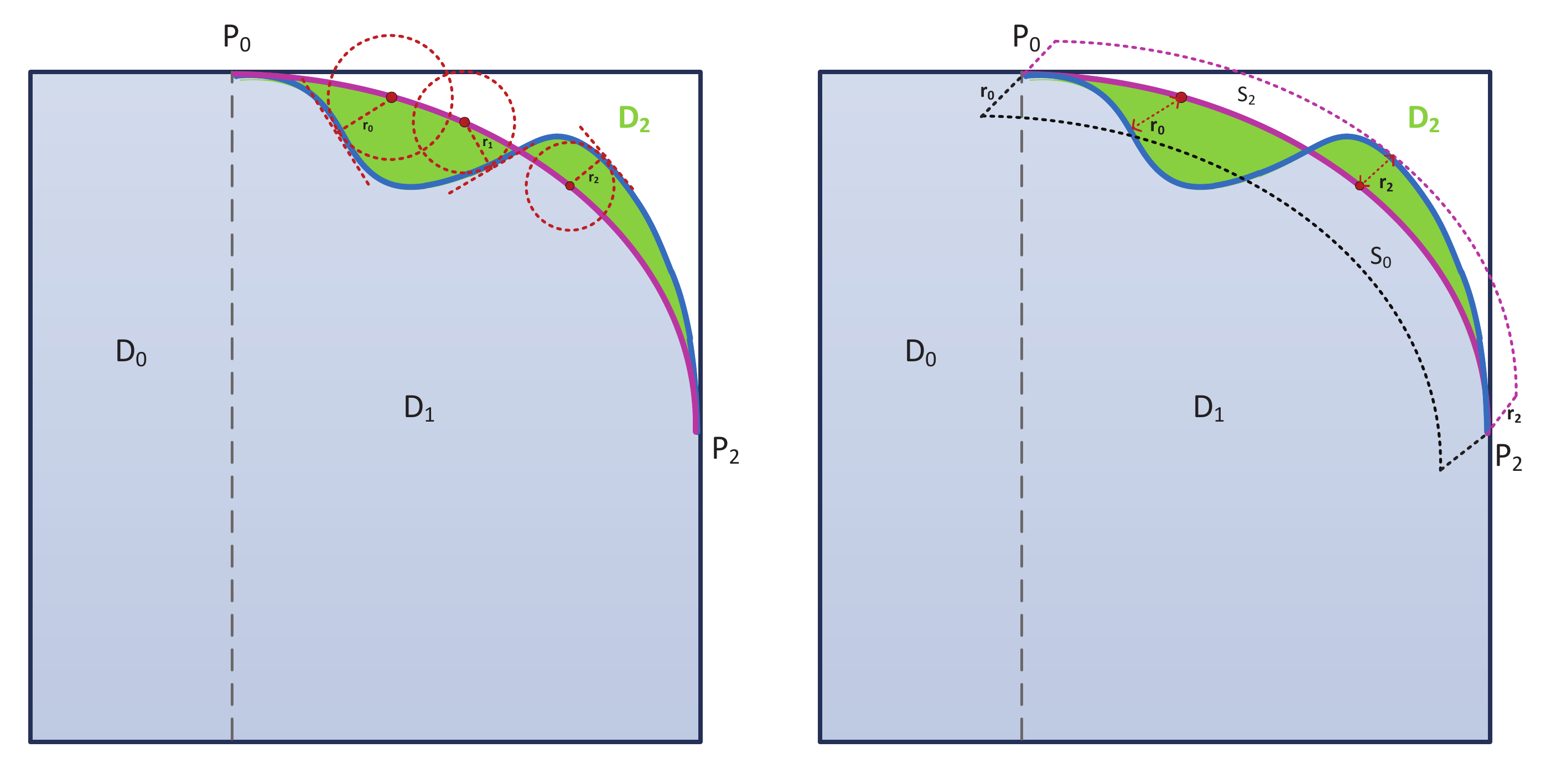}
	}
	\caption{Process of obtaining the area of $D_2$ (green). Left: Three sample points with their corresponding geometric distances $r_i,i=0,1,2$. Right: Area of $D_2$ is approximated by a narrow band. Its width is the longest geometric distance ($r_0$ on the left). The boundary of integration domain (blue), b\'{e}zier approximation to the boundary curve (magenta).
	}\label{fig:narrowband}
\end{figure}

To be more specific, let us take a close look at a specific boundary cell.
Figure \ref{fig:narrowband} (left) shows a boundary cell, where $P_0$ and $P_2$ are intersections of the domain to the cell grid. The domain in this cell is partitioned into three parts $D_0$, $D_1$ and $D_2$, where $D_0$ is a rectangle inside the domain, $D_1$ is a curved rectangle and $D_2$ is the region enclosed by the boundary curve and the approximated one. The former two regions will be considered during the numerical integration process. However, $D_2$ will be missed totally. This process leads to errors not only caused by numerical integration schemes, but also by the approximation.  
Here, the integration over $D_2$ is considered as the error guiding the subdivision.

According to the mean value theorem for integral, the error estimate over $D_2$ can be formulated as the following:
$$E_2=\iint_{D_2}|F|d\Omega\leq{M}\cdot S_{D_2},$$ where $M$ is the upper bound of absolute value of the integrand function $F$, and $S_{D_2}$ is the area of $D_2$. Since $F$ is continuous, the upper bound $M$ can always be obtained. Hence, our task turns to calculating the area of $D_2$.

A narrow band is used to give a pessimistic estimate on the area of $S_{D_2}$. As in Figure \ref{fig:narrowband} (right), the arc length of the B\'{e}zier curve in this cell is the length of the narrow band. To obtain the width of the band, the geometric distance between the boundary curve and the B\'{e}zier curve at several sample points are calculated, and the largest distance is chosen as the approximated width.  
To facilitate computation of width, Sampson distance is employed, which is the first order approximation to the geometric distance.

When the error estimate satisfies the scaled tolerance:
\begin{equation}\label{judgement}
E_2<\omega\tau
\end{equation}
then no more subdivisions are needed. The scale $\omega$ is the percentage of the area of the boundary cell occupies in the bounding box of the whole integration domain.

\begin{algorithm}[!h]
	\begin{algorithmic}[1]
		\caption{Get the numerical integration over implicitly defined domain $\int_{\Omega\cap{C}}{F}$}\label{lala2}
		\REQUIRE Domain $\Omega$ enclosing in $C=[a,b]\times[c,d]$, the smooth integrand $F$;\\		
		Tolerance $\tau$ (to guide the accuracy of the integral).
		\ENSURE The numerical result of $s=\int_\Omega{F}$
		
		\STATE $C$ is the  cells in the support set of the basis functions associated with $Pt$
		\STATE Initialization: let cell $C=[a,b]\times[c,d]$, $h =\max\{b-a,d-c\}$, we want to derive.  
		\STATE Initialization: set $h=b-a$, $g=d-c$, $s=0$.  
		\WHILE{$\max\{h,g\} \geq \epsilon$} 
		\STATE Classification of cell type (using interval arithmetic)
		
		\IF{ $C$ is an interior cell}
		\STATE Computing the integral over $C$, add the value to $s$;
		\ELSIF{ $C$ is a boundary cell}
		\STATE Using quadratic B\'{e}zier approximation to the implicit boundary, 
		\STATE Derive the local error estimate $E_2$ over this cell.
		\IF{ error criteria \eqref{judgement} is satisfied}
		\STATE Computing the integral with isoparametric integration, add it to $s$;
		\ELSE
		\STATE Subdivide $C$ into equal four subcells $C_i,i=1,2,3,4$; 
		\STATE Go to step $1$ calculate $\int_{\Omega\cap{C_i}}{F}$, add the four values to $s$.
		\ENDIF
		\ENDIF
		\ENDWHILE
		\STATE \textbf{end}
		\STATE return $s$.		
	\end{algorithmic}
\end{algorithm}

To conclude the method section, we here illustrate the main routine of the proposed method represented by Algorithm~\ref{lala2}.
The process begins with the overall bounding box of the full domain, clarifies the cell type. 
Once the specific type is determined, 
different types of cells lead to separate subroutines, following the subsequent approximations and error estimating criteria.
In order to improve the efficiency of the proposed method, one may combine corner classification together with interval classification to improve the efficiency.

\section{Numerical Experiments} \label{tests}

In this section, we will investigate the accuracy and performance of the proposed method, which uses interval arithmetic, the quadratic B\'{e}zier approximation to the boundary curve and a geometry-based local error estimate to guide the subdivision.

Four different domains are tested: annulus, a domain defined by splines, a cardioid domain and a cassini oval domain.
Among them, the first two consider domains with smooth boundaries, while the third concerns the domain with a singular point at its boundary, and the fourth examines a self-intersected domain.
We will test integrals with and without interval arithmetic to see the effect of interval arithmetic in preserving the correct topology of the domain of integration.

Comparisons are taken with another two methods. Both the two methods are hierarchy-based, with the prescribed subdivision length $1.2\cdot h/2^k$, where $h$ is the width of the bounding box of the domain, and $k$ shows the level of subdivision. 
Since linear and quadratic B\'{e}zier curves are used to approximate the boundary curve, respectively. Let us label the two methods separately as L and Q. In our method, the guiding tolerance $\tau=0.1^m, m\in N^+$. Accuracy of each method is quantified by the absolute error with the exact value, i.e., $|I-I_h|$, where $I_h$ is the numerical result and $I$ is the exact value obtained through Mathematica. 

\begin{exmp}
	Firstly, we will test the area of an annulus. The boundary of the annulus is defined by the zero level set of the implicit function $f(x,y)=0.04-(\sqrt{x^2+y^2}-0.6)^2$ on $[-1,1]^2$. The integrand is a unity function and the exact value of the integral is $\frac{12\pi}{25}$.

	The integration domain is highlighted in Figure \ref{fig:annulus1} (a). It can be seen clearly that, all the function values at corner points after subdividing the domain once are smaller than zero. This indicates that the subcells are exterior cells and hence the integration result is zero. However, with the help of interval arithmetic, the four subdomains are all boundary cells and the correct topology is thus attained.

	A group of comparison results (with the absolute error at magnitude $1E-5$) are depicted in Figure \ref{fig:annulus1} (b-d). It can be seen that both Q and our method use relatively smaller numbers of integration points compared with L method. In addition, since our method is error guided subdivision, the distribution of integration points is more adaptive compared with Q method (from the part enclosing the interior circle). 
	\begin{figure}[!h]
		\centerline{	
			\includegraphics[width=6in]{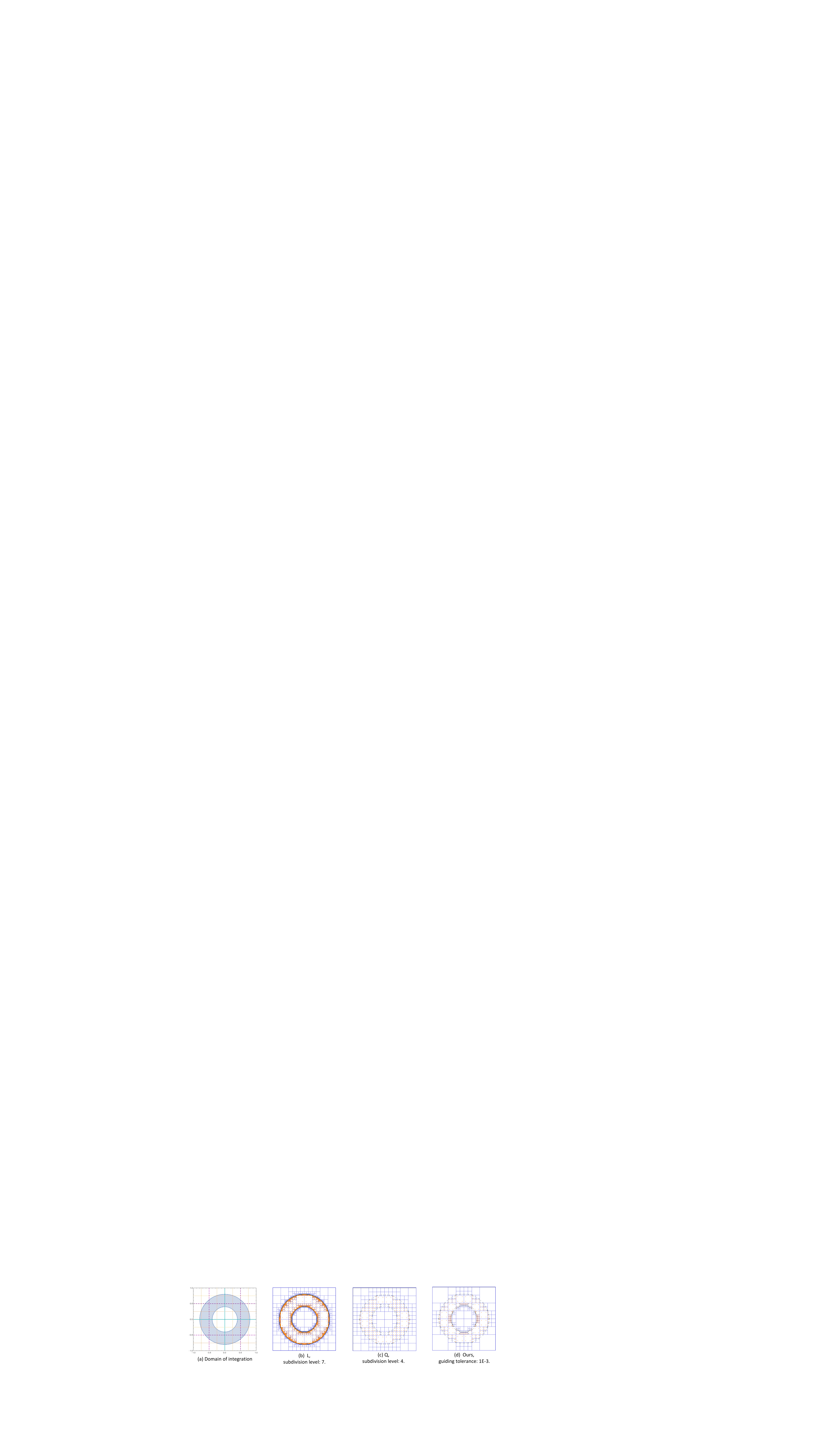}
		}\caption{Annulus area test. (a): the integration domain (blue). Lines in different colors are plotted to show the signs of function values at different levels. (b-d): the cell partitions and integration points of L, Q and our method, respectively. 
			The integration points in interior (boundary) cells are marked in orange (blue).}	\label{fig:annulus1}
	\end{figure}
	
		\begin{figure}[!h]
		\centerline{	
			\includegraphics[width=4in]{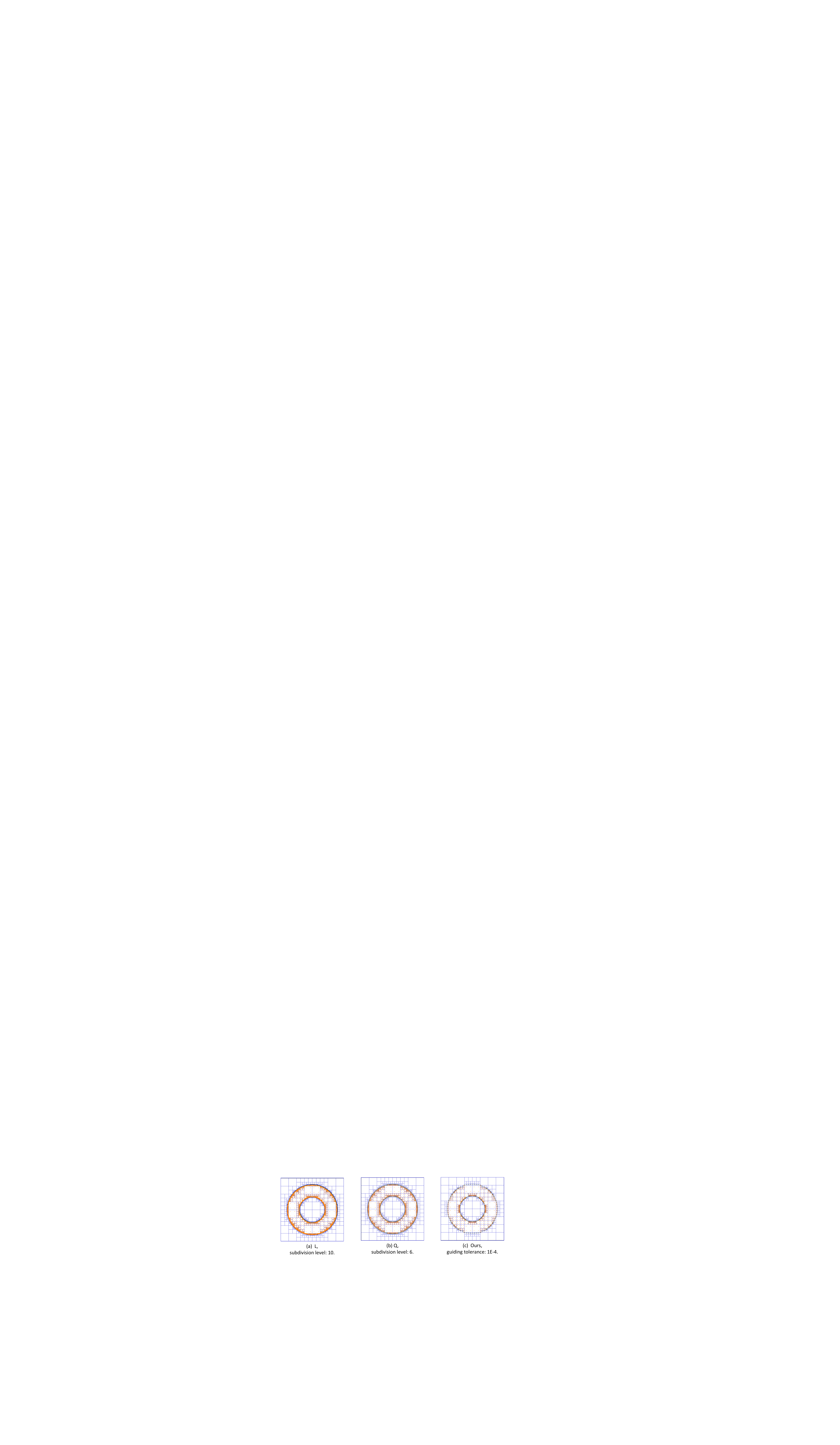}
		}\caption{Integratoin over annulus with integrand: $x^3y-xy+2.5$. 
			(a-c): the cell partitions and integration points of L, Q and our method, respectively. 
		}	\label{fig:annulus2}
	\end{figure}

	Now consider the integration over the annulus domain with integrand function $F(x,y)=x^3y-xy+2.5$. With two Gaussian points chosen in each direction, integrations over the interior cells are still accurate. In this way, the integration error is caused by the approximation to the boundary and the numerical quadrature over the boundary cells. Comparisons of different methods with absolute error at level $1E-7$ are shown in Figure \ref{fig:annulus2}.	From the figure we can see that our method uses the least numbers of integration points calculated among the three methods.

	Besides, the results of our method always satisfy the given guiding tolerances, as shown in Table \ref{table:annulus2}. That is, the pessimistic error estimate of our method leads to more accurate results (with absolute error smaller than the guiding tolerances), although there is no theoretical guarantee on the accuracy.	This is different from L and Q, where the smallest length is given to guide the uniform subdivision. Owing to this fact, we do not know how accurate the numerical integration derives. In the sense that we know how far the numerical method achieves, our method is better than the other two methods.

	As to the efficiency, here, we introduce a ratio to examine whether the geometry-based subdivision rule in our method works and the extent. Since the judgment (\ref{judgement}) is either true (no more subdivisions) or false (take another subdivision), we count the numbers of both cases of the condition (\ref{judgement}). Then let Cr be the ratio of the numbers of the satisfied error condition takes up in the whole judging process. In this way, the obtained Cr measures the calculations saved by the subdivision rule (the satisfied error condition do not need more levels of subdivisions).
	
		\begin{table}[!h]\small
		\centering
		\caption{Integration over annulus with different integrands. Cr is the ratio 
			of the number of times when the condition (\ref{judgement}) is satisfied to the number of times of judgments on the condition (\ref{judgement}).}\label{table:annulus2}
		\begin{tabular}{|c|ccccc|ccccc|}
			\hline
			\multirow{2}{*}{Tolerance}	&\multicolumn{5}{c|}{Integrand: 1}&\multicolumn{5}{c|}{Integrand: $x^3y-xy+2.5$}\\ 
			\cline{2-11}	
			& Error&Time(ms)&\# In&\# Bd&Cr&Error&Time(ms)&\#In&\# Bd&Cr\\ \hline
			$1.00E-01$&	$2.31E-03$&	$0.2$&	$0$&	$24$&	$0.86$ &	$5.78E-03$&	$0.25$&	$0$&	$24$&	$0.86$ \\ \hline
			$1.00E-02$&	$1.10E-04$&	$0.45$&	$8$&	$56$&	$0.74$ &	$2.75E-04$&	$0.55$&	$8$&	$56$&	$0.74$ \\ \hline
			$1.00E-03$&	$1.17E-05$&	$0.8$&	$32$&	$104$&	$0.67$ &	$2.89E-05$&	$1$&	$32$&	$104$&	$0.67$ \\ \hline
			$1.00E-04$&	$5.18E-06$&	$1.45$&	$72$&	$168$&	$0.63$ &	$4.25E-07$&	$2.7$&	$132$&	$264$&	$0.61$ \\ \hline
			$1.00E-05$&	$3.83E-08$&	$3.3$&	$196$&	$376$&	$0.60$ &	$2.73E-07$&	$4.15$&	$220$&	$472$&	$0.61$ \\ \hline
			$1.00E-06$&	$5.81E-09$&	$5.85$&	$380$&	$728$&	$0.60$ &	$3.32E-07$&	$8.55$&	$508$&	$1048$&	$0.61$ \\ \hline
		\end{tabular}
	\end{table}

	It is seen from Table \ref{table:annulus2} that, all the values of Cr are not less than $0.6$. This indicates that the error-guided subdivision scheme saves calculation of more than half percent of cells compared with Q to the same level of subdivision. 
	
\end{exmp}

\begin{exmp}
	Next, we explore the integration over a domain defined by a bi-quadratic tensor product B-spline function with two different integrands. The knot sequence in each direction is $\{0,0,0,0.3, 0.7, 1.0,1.0,1.0\}$, from which we get an initial partition of the domain, as depicted in Figure \ref{fig:rabbitiacmp} (a) and (b). The function is represented as $$f(x,y)=\sum\limits_{i=0}^{4}\sum\limits_{j=0}^{4}{P_{ij}N_i^3(x)N_j^3(y)},$$ where $N_k^3(t)$ is B-spline function, and the coefficients matrix is 
$$
P=(P_{ij})_{5\times5}=\left[\begin{matrix}
-1& -1.5& 1& -8&-4\\
-1& 2& 1& 3.2& -1\\
-2& 3& -2.2& 2.5& -1\\
-1& -1& 2& 3& 0.8\\
-0.2& 1.5& -1& 0.8& 0.3  
\end{matrix}
\right]. $$

Let us firstly set the integrand to be 1 and consider the area of the domain. Figure \ref{fig:rabbitiacmp} (c-e) show the cell partitions and integration points of L without using interval arithmetic. It is easily observed that two boundary cells (in the middle and in the middle right of the domain) are misjudged as interior cells thus no more subdivisions are taken. (f) is the result of L with interval arithmetic (IA) corresponding to (d). Clearly, we see that the integration domain is kept complete. Furthermore, (f-h) depict the results of L, Q and our method with absolute error at $1E-4$ using interval arithmetic. It is observed that all the three methods achieve correct topology. What is more, our method uses the least numbers of integration points. 

\begin{figure}[!h]
	\centerline{	
		\includegraphics[width=5.3in]{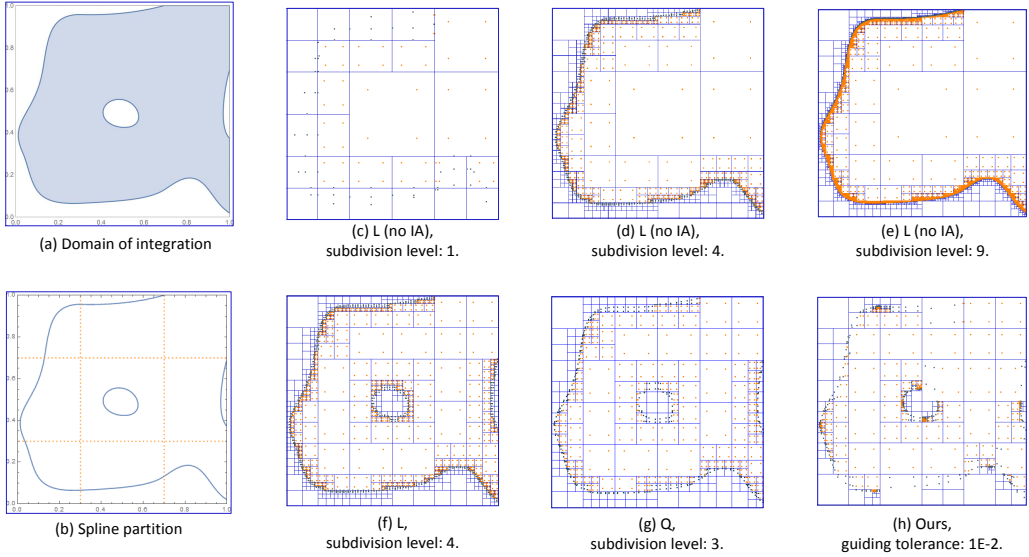}
	}		\caption{Comparisons of integral of the area of the spline domain.
		}\label{fig:rabbitiacmp}	
\end{figure}

Detailed comparisons for integrations over the spline domain with and without interval arithmetic (IA) are listed in Table \ref{table:rabbit0} and Table \ref{table:rabbit1}, respectively. Different integrands are tested, and each case derives the consistent results. In the former table, all the results have large errors (magnitude $1E-2$) due to the incorrect boundary detection; while in the latter table, since interval arithmetic is used, errors will decrease.

\begin{table}[!h]\small
	\centering
	\caption{Comparisons of integration over the spline domain with different integrands (no IA). }\label{table:rabbit0}
	\begin{tabular}{|c|c|c|c|c|}
		\hline
		\multirow{2}{*}{Method}	&\multicolumn{2}{c|}{Integrand: 1}&\multicolumn{2}{c|}{Integrand: $x^3y-xy+2.5$}\\ 
		\cline{2-5}	
	& Error	& Error	& Error	& Error\\ \hline
		L&$2.39E-02$&$2.40E-02$&$5.68E-02$&$5.68E-02$\\
		Q&$2.40E-02$&$2.40E-02$&$5.68E-02$&$5.68E-02$\\
		Our method&$2.39E-02$&$2.40E-02$&$5.67E-02$&$5.68E-02$ \\ \hline
	\end{tabular}
\end{table} 

\begin{table}[!h]\small
	\centering
	\caption{Comparisons of integration over the spline domain with different integrands.}\label{table:rabbit1}
	\begin{tabular}{|c|cccc|cccc|}
		\hline
		\multirow{2}{*}{Method}	&\multicolumn{4}{c|}{Integrand: 1}&\multicolumn{4}{c|}{Integrand: $x^3y-xy+2.5$}\\ 
		\cline{2-9}	
			& Error&Time(ms)&\# In&\# Bd&Error&Time(ms)&\#In&\# Bd\\ \hline
		L&  	$2.74E-04$& 	$8.75$& 	$175$& 	$180$&  	$1.18E-05$& 	$67$& 	$1458$& 	$1429$\\ 
		Q& 	$6.34E-04$& 	$10.5$& 	$93$& 	$87$& 	$1.53E-05$& 	$43$& 	$359$& 	$360$\\
		Our method& 	$1.39E-04$& 	$16.6$& 	$78$& 	$98$& 	$2.51E-05$& 	$22$& 	$108$& 	$132$\\ \hline
		L& 		$2.41E-06$& 	$133.05$& 	$2962$& 	$2856$&  	$4.92E-06$& 	$131$& 	$2962$& 	$2856$\\ 
		Q& 		$2.09E-06$& 	$83.7$& 	$721$& 	$718$&  	$4.21E-06$& 	$85$& 	$721$& 	$718$\\ 
		Our method&  	$3.21E-06$& 	$28.05$& 	$126$& 	$168$&  	$1.98E-06$& 	$42$& 	$181$& 	$235$\\ \hline 	\end{tabular}
\end{table}

From Table \ref{table:rabbit1}, we also find that, among the three methods for most cases, L is the most time consuming and calculates the highest numbers of cells, and that of Q is smaller than L. However, our method is always the least time-consuming and saves the calculation numbers of cells most.

	\begin{table}[!h]\small
	\centering
	\caption{Integration over the spline domain with different integrands and tolerances.}\label{table:rabbit2}
	\begin{tabular}{|c|ccccc|ccccc|}
		\hline
		\multirow{2}{*}{Tolerance}	&\multicolumn{5}{c|}{Integrand: 1}&\multicolumn{5}{c|}{Integrand: $x^3y-xy+2.5$}\\ 
		\cline{2-11}	
		& Error&Time(ms)&\# In&\# Bd&Cr&Error&Time(ms)&\#In&\# Bd&Cr\\ \hline
		$1.00E-01$&	$4.25E-04$&	$8.45$&	$38$&	$52$&	$0.64$ &	$1.97E-03$&	$12$&	$54$&	$74$&	$0.62$ \\ \hline
		$1.00E-02$&	$1.39E-04$&	$16.6$&	$78$&	$98$&	$0.59$ &	$2.51E-05$&	$23$&	$108$&	$132$&	$0.58$ \\ \hline
		$1.00E-03$&	$3.21E-06$&	$28.05$&	$126$&	$168$&	$0.58$ &	$1.98E-06$&	$41$&	$181$&	$235$&	$0.57$ \\ \hline
		$1.00E-04$&	$9.87E-07$&	$56.2$&	$258$&	$337$&	$0.57$ &	$2.07E-06$&	$76$&	$337$&	$456$&	$0.57$ \\ \hline
		$1.00E-05$&	$1.32E-06$&	$109.25$&	$462$&	$669$&	$0.57$&	$2.43E-06$&	$142$&	$573$&	$871$&	$0.58$ \\ \hline
	\end{tabular}
\end{table}

Table \ref{table:rabbit2} tabulates the results of our method under different guiding tolerances. From the table, we can see that the errors of our method in all tests are smaller than the given tolerances. What is more, all the values of Cr for different cases presented in the table are larger than $0.56$, showing that the subdivision rule saves more than half numbers of subdivisions. It thus demonstrates the efficiency of our local error estimate.
\end{exmp}

\begin{exmp}
We now examine the numerical integration with two integrands over the cardioid domain. There is a singular point at the boundary curve, as shown in Figure \ref{fig:cardcmp} (a). With the background lines, it is easy observed that the cell partition just using the corner points, even together with center points is not enough to get the correct shape.

\begin{figure}[!h]
	\centerline{	
		\includegraphics[width=5.8in]{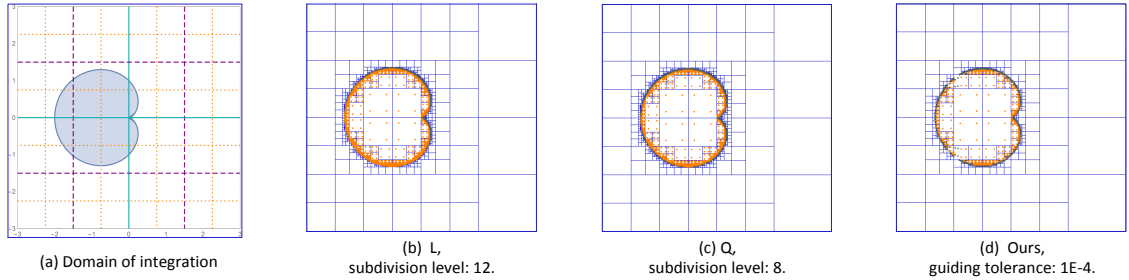}
	}		\caption{Integration over cardioid using different approximations to the boundary with interval arithmetic.}
	\label{fig:cardcmp}	
\end{figure}

Figure \ref{fig:cardcmp} (b-d) illustrate the results of absolute errors at $1E-6$ using different methods. Since there are two branches at the singular point $(0,0)$ and both tangents at the point are zero, we just split the domain into four parts, and handle each case separately. As we can see, the number of integration points of our method is the smallest.

\begin{table}[!h]\small
	\centering
	\caption{Integration over the cardioid domain with different integrands and tolerances.}\label{table:card1}
	\begin{tabular}{|c|ccccc|ccccc|}
		\hline
		\multirow{2}{*}{Tolerance}	&\multicolumn{5}{c|}{Integrand: 1}&\multicolumn{5}{c|}{Integrand: $x^3y-xy+3$}\\ 
		\cline{2-11}	
		& Error&Time(ms)&\# In&\# Bd&Cr&Error&Time(ms)&\#In&\# Bd&Cr\\ \hline
		$1.00E-01$ &	$1.47E-02$&	$0.1$&	$12$&	$22$&	$0.58$ &	$9.94E-03$&	$0.1$&	$22$&	$32$&	$0.57$ \\ \hline
		$1.00E-02$ &	$2.48E-04$&	$0.1$&	$40$&	$50$&	$0.58$ &	$6.98E-04$&	$0.2$&	$56$&	$60$&	$0.57$ \\ \hline
		$1.00E-03$&	$1.09E-04$&	$0.3$&	$80$&	$72$&	$0.47$ &	$3.50E-05$&	$0.5$&	$110$&	$126$&	$0.50$ \\ \hline
		$1.00E-04$&	$3.05E-06$&	$0.9$&	$192$&	$182$&	$0.48$ &	$2.65E-06$&	$1.4$&	$302$&	$264$&	$0.45$ \\ \hline
		$1.00E-05$&	$8.48E-08$&	$1.7$&	$404$&	$380$&	$0.45$ &	$1.71E-07$&	$2.7$&	$642$&	$570$&	$0.42$ \\ \hline
		$1.00E-06$&	$2.05E-08$&	$3.5$&	$916$&	$754$&	$0.36$ &	$1.99E-08$&	$5.9$&	$1512$&	$1164$&	$0.32$ \\ \hline
	\end{tabular}
\end{table} 

The detailed results of our method under different tolerances are presented in Table \ref{table:card1}. It can be seen that our method yields errors smaller than the required tolerances in all tests.

As to the efficiency, different from the cases of the above two smooth domains, the results here have smaller Cr. The reason is that the geometry may vary much near the singular points, which causes more levels of subdivisions to satisfy the local error condition. Since Cr is larger than 0.3, it indicates that although the ratio is not that large, it still saves calculations for the singular case. 
\end{exmp}

\begin{exmp}
	In this example, we consider the numerical integration over a cassini oval domain with self-intersection boundary, as depicted in Figure \ref{fig:eightcmp} (a). The function defining the domain is $f(x,y)=0.98(x^2 - y^2) - (x^2 + y^2)^2$. The integrands are the same with the above cardioid test. From the blue background lines in the figure, it can be seen that all the function values of the corner points of the four subcells are not larger than zero. Therefore, if not using interval arithmetic, all the subcells are indicated to be exterior subcells. Hence, we could only get integration value zero, no matter what the integrand is. However, with interval arithmetic, the four subcells are recognized as boundary cells, which suggests more levels of calculations. In this way, the correct topology is preserved and accurate integration results are derived. 
 
 		\begin{figure}[!h]
			\centerline{	
				\includegraphics[width=5.8in]{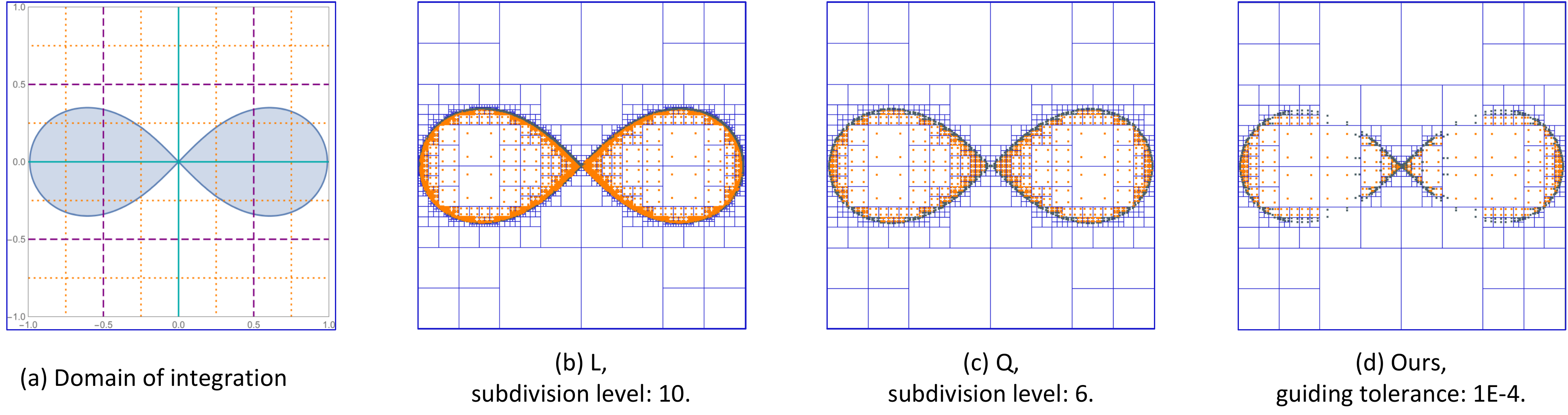}
			}		\caption{Integration over cassini oval using different methods with interval arithmetic.}
			\label{fig:eightcmp}	
		\end{figure}
	
	Figure \ref{fig:eightcmp} (b-d) show the results of area test with absolute errors at $1E-6$ using different methods. Following the same fashion as the former example, all the three methods lead to accurate results, and our method uses the least number of integration points.
	
	Table \ref{table:eight} tabulates the comparison of results with different integrands and different methods with interval arithmetic. The results show that our method behaves better than L and Q. That is, for the same integrand, to the same grade of absolute error, our method usually uses the smallest number of integration points, well costs the shortest time.
	
\begin{table}[!h]\small
	\centering
	\caption{Comparisons of integration over the cassini oval domain with different integrands and methods.}\label{table:eight}
	\begin{tabular}{|c|cccc|cccc|}
		\hline
		\multirow{2}{*}{Method}	&\multicolumn{4}{c|}{Integrand: 1}&\multicolumn{4}{c|}{Integrand: $x^3y-xy+3$}\\ 
		\cline{2-9}	
		& Error&Time(ms)&\# In&\# Bd&Error&Time(ms)&\#In&\# Bd\\ \hline
	L	&	$4.55E-05$	&	$1.55$	&	$840$	&	$860$	&	$1.37E-04$	&	$1.55$	&	$840$	&	$860$	\\
	Q	&	$5.46E-05$	&	$0.3$	&	$76$	&	$104$	&	$1.64E-04$	&	$0.3$	&	$76$	&	$104$	\\
	Our method	&	$3.11E-05$	&	$0.25$	&	$20$	&	$92$	&	$1.09E-04$	&	$0.3$	&	$12$	&	$72$	\\  \hline 
	L	&	$2.64E-06$	&	$6.65$	&	$3568$	&	$3460$	&	$7.93E-06$	&	$6.1$	&	$3568$	&	$3460$	\\
	Q	&	$3.53E-06$	&	$0.6$	&	$184$	&	$216$	&	$1.20E-06$	&	$1.25$	&	$408$	&	$432$	\\
	Our method	&	$1.15E-06$	&	$0.55$	&	$52$	&	$180$	&	$6.45E-06$	&	$0.45$	&	$40$	&	$144$	\\
	 \hline 
	\end{tabular}
\end{table}
\end{exmp}

\section{Conclusion} \label{con}
We have developed a method for numerical integration over implicitly defined domains. In this method, interval arithmetic is adopted under the hierarchical frame, together with a geometric based local error estimate to guide the subdivision. The results show that this method preserves the correct topology, and leads to accurate and efficient results for continuous integrands. We believe that this work gives one way for accurate integration that considering the geometry information.

There are several directions for further research. The extension of this idea to 3D integration is straightforward. Considering that affine arithmetic is efficient when detecting the boundary of an implicit domain, an interesting future topic would be to adopt it when performing numerical quadrature. Besides, more careful error estimate could be taken to further improve the accuracy and efficiency. Applications in engineering would also be one of the potential directions.

\section*{Acknowledgment}
We would like to thank the anonymous reviewers and our laboratory group for helpful discussions and comments. The
work is supported by the National Natural Science Foundation of China (No.11771420).

\bibliographystyle{unsrt} 
\bibliography{m2dquadrature} 

\begin{thebibliography}{10}

\bibitem{huang2013intersection}
Pu~Huang, Charlie~CL Wang, and Yong Chen.
\newblock Intersection-free and topologically faithful slicing of implicit
  solid.
\newblock {\em Journal of Computing and Information Science in Engineering},
  13(2):021009, 2013.

\bibitem{wegst2015bioinspired}
Ulrike~GK Wegst, Hao Bai, Eduardo Saiz, Antoni~P Tomsia, and Robert~O Ritchie.
\newblock Bioinspired structural materials.
\newblock {\em Nature materials}, 14(1):23, 2015.

\bibitem{cottrell2009isogeometric}
J~Austin Cottrell, Thomas~JR Hughes, and Yuri Bazilevs.
\newblock {\em Isogeometric analysis: toward integration of CAD and FEA}.
\newblock John Wiley \& Sons, 2009.

\bibitem{upreti2017signed}
K~Upreti and G~Subbarayan.
\newblock Signed algebraic level sets on nurbs surfaces and implicit boolean
  compositions for isogeometric cad--cae integration.
\newblock {\em Computer-Aided Design}, 82:112--126, 2017.

\bibitem{dokken2018trivariate}
Tor Dokken, Vibeke Skytt, and Oliver Barrowclough.
\newblock Trivariate spline representations for computer aided design and
  additive manufacturing.
\newblock {\em arXiv preprint arXiv:1803.05756}, 2018.

\bibitem{hollig2003finite}
Klaus Hollig.
\newblock {\em Finite element methods with B-splines}, volume~26.
\newblock Siam, 2003.

\bibitem{rank2012geometric}
Ernst Rank, Martin Ruess, Stefan Kollmannsberger, Dominik Schillinger, and
  Alexander D{\"u}ster.
\newblock Geometric modeling, isogeometric analysis and the finite cell method.
\newblock {\em Computer Methods in Applied Mechanics and Engineering},
  249:104--115, 2012.

\bibitem{xu2017improved}
Jinlan Xu, Ningning Sun, Laixin Shu, Timon Rabczuk, and Gang Xu.
\newblock An improved isogeometric analysis method for trimmed geometries.
\newblock {\em arXiv preprint arXiv:1707.00323}, 2017.

\bibitem{bartovn2017gauss}
Michael Barto{\v{n}} and Victor~Manuel Calo.
\newblock Gauss--galerkin quadrature rules for quadratic and cubic spline
  spaces and their application to isogeometric analysis.
\newblock {\em Computer-Aided Design}, 82:57--67, 2017.

\bibitem{barendrecht2018efficient}
Pieter~J Barendrecht, Michael Barto{\v{n}}, and Ji{\v{r}}{\'\i} Kosinka.
\newblock Efficient quadrature rules for subdivision surfaces in isogeometric
  analysis.
\newblock {\em Computer Methods in Applied Mechanics and Engineering}, 2018.

\bibitem{xu2017isogeometric}
Gang Xu, Tsz-Ho Kwok, and Charlie~CL Wang.
\newblock Isogeometric computation reuse method for complex objects with
  topology-consistent volumetric parameterization.
\newblock {\em Computer-Aided Design}, 91:1--13, 2017.

\bibitem{dolbow1999numerical}
John Dolbow and Ted Belytschko.
\newblock Numerical integration of the galerkin weak form in meshfree methods.
\newblock {\em Computational mechanics}, 23(3):219--230, 1999.

\bibitem{peskin2002immersed}
Charles~S Peskin.
\newblock The immersed boundary method.
\newblock {\em Acta numerica}, 11:479--517, 2002.

\bibitem{li2006immersed}
Zhilin Li and Kazufumi Ito.
\newblock {\em The immersed interface method: numerical solutions of PDEs
  involving interfaces and irregular domains}, volume~33.
\newblock Siam, 2006.

\bibitem{moes1999finite}
Nicolas Mo{\"e}s, John Dolbow, and Ted Belytschko.
\newblock A finite element method for crack growth without remeshing.
\newblock {\em International journal for numerical methods in engineering},
  46(1):131--150, 1999.

\bibitem{sukumar2001modeling}
Natarajan Sukumar, David~L Chopp, Nicolas Mo{\"e}s, and Ted Belytschko.
\newblock Modeling holes and inclusions by level sets in the extended
  finite-element method.
\newblock {\em Computer methods in applied mechanics and engineering},
  190(46-47):6183--6200, 2001.

\bibitem{cheng2010higher}
Kwok~Wah Cheng and Thomas-Peter Fries.
\newblock Higher-order xfem for curved strong and weak discontinuities.
\newblock {\em International Journal for Numerical Methods in Engineering},
  82(5):564--590, 2010.

\bibitem{Edalat1999Numerical}
Abbas Edalat and Marko Krznaric.
\newblock {\em Numerical Integration with Exact Real Arithmetic}.
\newblock Springer Berlin Heidelberg, 1999.

\bibitem{berntsen1991adaptive}
Jarle Berntsen, Terje~O Espelid, and Alan Genz.
\newblock An adaptive algorithm for the approximate calculation of multiple
  integrals.
\newblock {\em ACM Transactions on Mathematical Software (TOMS)},
  17(4):437--451, 1991.

\bibitem{davis2007methods}
Philip~J Davis and Philip Rabinowitz.
\newblock {\em Methods of numerical integration}.
\newblock Courier Corporation, 2007.

\bibitem{muller2013highly}
Bj{\"o}rn M{\"u}ller, Florian Kummer, and Martin Oberlack.
\newblock Highly accurate surface and volume integration on implicit domains by
  means of moment-fitting.
\newblock {\em International Journal for Numerical Methods in Engineering},
  96(8):512--528, 2013.

\bibitem{press2007numerical}
William~H Press.
\newblock {\em Numerical recipes 3rd edition: The art of scientific computing}.
\newblock Cambridge university press, 2007.

\bibitem{shapiro2002architecture}
Vadim Shapiro and I~Tsukanov.
\newblock The architecture of sage--a meshfree system based on rfm.
\newblock {\em Engineering with Computers}, 18(4):295--311, 2002.

\bibitem{tornberg2004numerical}
Anna-Karin Tornberg and Bj{\"o}rn Engquist.
\newblock Numerical approximations of singular source terms in differential
  equations.
\newblock {\em Journal of Computational Physics}, 200(2):462--488, 2004.

\bibitem{rvachev1994numerical}
VL~Rvachev, AN~Shevchenko, and VV~Veretel'nik.
\newblock Numerical integration software for projection and projection-grid
  methods.
\newblock {\em Cybernetics and Systems Analysis}, 30(1):154--158, 1994.

\bibitem{thiagarajan2014adaptively}
Vaidyanathan Thiagarajan and Vadim Shapiro.
\newblock Adaptively weighted numerical integration over arbitrary domains.
\newblock {\em Computers \& Mathematics with Applications}, 67(9):1682--1702,
  2014.

\bibitem{thiagarajan2017shape}
Vaidyanathan Thiagarajan.
\newblock {\em Shape aware quadratures}.
\newblock The University of Wisconsin-Madison, 2017.

\bibitem{thiagarajan2016adaptively}
Vaidyanathan Thiagarajan and Vadim Shapiro.
\newblock Adaptively weighted numerical integration in the finite cell method.
\newblock {\em Computer Methods in Applied Mechanics and Engineering},
  311:250--279, 2016.

\bibitem{olshanskii2016numerical}
Maxim~A Olshanskii and Danil Safin.
\newblock Numerical integration over implicitly defined domains for higher
  order unfitted finite element methods.
\newblock {\em Lobachevskii Journal of Mathematics}, 37(5):582--596, 2016.

\bibitem{engwer2016geometric}
Christian Engwer and Andreas N{\"u}{\ss}ing.
\newblock Geometric integration over irregular domains with topologic
  guarantees.
\newblock {\em arXiv preprint arXiv:1601.03597}, 2016.

\bibitem{hollig2015programming}
Klaus H{\"o}llig and J{\"o}rg H{\"o}rner.
\newblock Programming finite element methods with weighted b-splines.
\newblock {\em Computers \& Mathematics with Applications}, 70(7):1441--1456,
  2015.

\bibitem{thiagarajan2018shape}
Vaidyanathan Thiagarajan and Vadim Shapiro.
\newblock Shape aware quadratures.
\newblock {\em Journal of Computational Physics}, 2018.

\bibitem{saye2015high}
RI~Saye.
\newblock High-order quadrature methods for implicitly defined surfaces and
  volumes in hyperrectangles.
\newblock {\em SIAM Journal on Scientific Computing}, 37(2):A993--A1019, 2015.

\bibitem{drescher2017high}
Lukas Drescher, Holger Heumann, and Kersten Schmidt.
\newblock A high order method for the approximation of integrals over
  implicitly defined hypersurfaces.
\newblock {\em SIAM Journal on Numerical Analysis}, 55(6):2592--2615, 2017.

\bibitem{gautschi1989gauss}
Walter Gautschi and Sotorios~E Notaris.
\newblock Gauss--kronrod quadrature formulae for weight functions of
  bernstein--szeg{\"o} type.
\newblock {\em Journal of Computational and Applied Mathematics},
  25(2):199--224, 1989.

\bibitem{moore1966interval}
R~Moore.
\newblock {\em Interval arithmetic}.
\newblock Prentice-Hall Englewood Cliffs, 1966.

\bibitem{mitchell1990robust}
Don~P Mitchell.
\newblock Robust ray intersection with interval arithmetic.
\newblock In {\em Proceedings of Graphics Interface}, volume~90, pages 68--74,
  1990.

\bibitem{moore2014reliability}
Ramon~E Moore.
\newblock {\em Reliability in computing: the role of interval methods in
  scientific computing}, volume~19.
\newblock Elsevier, 2014.

\bibitem{martin2002comparison}
Ralph Martin, Huahao Shou, Irina Voiculescu, Adrian Bowyer, and Guojin Wang.
\newblock Comparison of interval methods for plotting algebraic curves.
\newblock {\em Computer Aided Geometric Design}, 19(7):553--587, 2002.

\bibitem{shou2003modified}
Huahao Shou, Hongwei Lin, Ralph Martin, and Guojin Wang.
\newblock Modified affine arithmetic is more accurate than centered interval
  arithmetic or affine arithmetic.
\newblock In {\em Mathematics of Surfaces}, pages 355--365. Springer, 2003.

\bibitem{moore1979methods}
Ramon~E Moore.
\newblock {\em Methods and applications of interval analysis}, volume~2.
\newblock Siam, 1979.

\bibitem{gomes2009implicit}
Abel Gomes, Irina Voiculescu, Joaquim Jorge, Brian Wyvill, and Callum
  Galbraith.
\newblock {\em Implicit curves and surfaces: mathematics, data structures and
  algorithms}.
\newblock Springer Science \& Business Media, 2009.

\bibitem{caprani1981integration}
Ole Caprani, Kaj Madsen, and Louis~B Rall.
\newblock Integration of interval functions.
\newblock {\em SIAM Journal on Mathematical Analysis}, 12(3):321--341, 1981.

\bibitem{wolfe1998interval}
MA~Wolfe.
\newblock Interval enclosures for a certain class of multiple integrals.
\newblock {\em Applied mathematics and computation}, 96(2-3):145--159, 1998.

\bibitem{farin2002curves}
Gerald~E Farin.
\newblock {\em Curves and surfaces for CAGD: a practical guide}.
\newblock Morgan Kaufmann, 2002.

\bibitem{hartshorne2013algebraic}
Robin Hartshorne.
\newblock {\em Algebraic geometry}, volume~52.
\newblock Springer Science \& Business Media, 2013.

\end{thebibliography}

\end{document}